\documentclass[12pt]{article}

\newtheorem{theorem}{Theorem}[section]
\newtheorem{lemma}[theorem]{Lemma}
\newtheorem{proposition}[theorem]{Proposition}
\newtheorem{corollary}[theorem]{Corollary}

\newenvironment{poof}{\textit{Proof:  }}{
~\hfill\rule{2mm}{3mm}\vspace{.2in}}

\usepackage{amssymb,amsmath,tabularx,graphicx,float,placeins}
\usepackage{tikz}
\usetikzlibrary[calc,intersections,through,backgrounds,arrows,decorations.pathmorphing]

\def\Acal{\mathcal{A}}

\def\Pcal{\mathcal{P}}

\def\Tcal{\mathcal{T}}
\def\Wcal{\mathcal{W}}

\def\Rbb{\mathbb{R}}

\def\ra{\rightarrow}
\def\Ra{\Rightarrow}

\def\ov{\overline}

\def\pr{\prime}
\def\prpr{\prime\prime}

\DeclareMathOperator{\lk}{lk}
\DeclareMathOperator{\dl}{dl}
\DeclareMathOperator{\st}{st}

\DeclareMathOperator{\Con}{Con}
\DeclareMathOperator{\Cg}{Cg}

\begin{document}

\title{Crosscut-simplicial Lattices}
\date{}
\author{Thomas McConville}

\maketitle

\begin{abstract}
We call a finite lattice crosscut-simplicial if the crosscut complex of every atomic interval is equal to the boundary of a simplex.  Every interval of such a lattice is either contractible or homotopy equivalent to a sphere.  Recently, Hersh and M\'esz\'aros introduced SB-labellings and proved that if a lattice has an SB-labelling then it is crosscut-simplicial.  Some known examples of lattices with a natural SB-labelling include the join-distributive lattices, the weak order of a Coxeter group, and the Tamari lattice.  Generalizing these three examples, we prove that every meet-semidistributive lattice is crosscut-simplicial, though we do not know whether all such lattices admit an SB-labelling.  While not every crosscut-simplicial lattice is meet-semidistributive, we prove that these properties are equivalent for chamber posets of real hypeplane arrangements.
\end{abstract}

\section{Introduction}\label{sec_introduction}

Many familiar posets have a M\"obius function that only takes values in the set $\{1,-1,0\}$.  To explain this occurrence, Hersh and M\'esz\'aros introduced SB-labellings, a labelling of the covering relations of a lattice which ensures that every interval is either contractible or homotopy equivalent to a sphere \cite{hersh.meszaros:sb}.  As usual, the topology associated to a poset is that of its order complex, the simplicial complex of chains.  A general fact that we use repeatedly is that a poset is homotopy equivalent to any of its crosscut complexes, which we review in Section \ref{sec_CSL}.  Unless specified otherwise, we always take the crosscut of atoms to define the crosscut complex of a lattice.

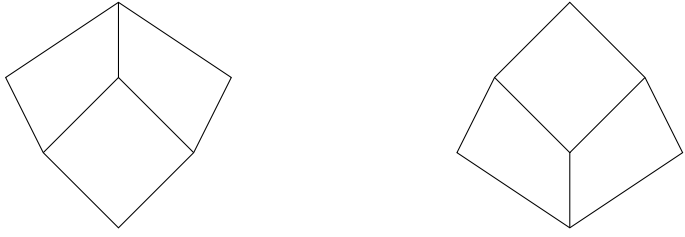
\begin{figure}
\begin{centering}
\begin{tikzpicture}
\begin{scope}
\draw (0,0) -- (-1,1) -- (-1.5,2) -- (0,3) -- (1.5,2) -- (1,1) -- cycle;
\draw (0,3) -- (0,2) -- (-1,1)
    (0,2) -- (1,1);
\end{scope}
\begin{scope}[xshift=6cm]
\draw (0,0) -- (-1.5,1) -- (-1,2) -- (0,3) -- (1,2) -- (1.5,1) -- cycle;
\draw (0,0) -- (0,1) -- (-1,2)
    (0,1) -- (1,2);
\end{scope}
\end{tikzpicture}
\caption{\label{fig_ex_CSL}\scriptsize The lattice on the left is crosscut-simplicial, while the lattice on the right is not.}
\end{centering}
\end{figure}

Hersh and M\'esz\'aros proved that lattices with SB-labellings are crosscut-simplicial.  We say a lattice is \emph{crosscut-simplicial} if for any interval $[x,y]$, the join of any proper subset of atoms of $[x,y]$ is not equal to $y$; see Figure \ref{fig_ex_CSL}.  Equivalently, a lattice is crosscut-simplicial if the crosscut complex on the atoms of any atomic interval is the boundary of a simplex.  In particular, a crosscut-simplicial lattice has every interval either contractible or homotopy equivalent to a sphere.

Two operations that preserve the crosscut-simplicial property are lattice quotients and doublings at order filters.  Lattice quotients and doublings are recalled in Sections \ref{sec_quotient} and \ref{sec_doubling}.

\begin{theorem}\label{thm_main_quotient}
Let $L$ be a crosscut-simplicial lattice.
\begin{enumerate}
\item If $\Theta$ is a lattice congruence of $L$, then the quotient $L/\Theta$ is crosscut-simplicial.
\item If $C$ is an order filter of $L$, then the doubling of $L$ at $C$ is crosscut-simplicial.
\end{enumerate}
\end{theorem}

The weak order of a finite Coxeter system inherits an SB-labelling from its Cayley graph (\cite{hersh.meszaros:sb} Theorem 5.3).  The weak order was originally proved to be crosscut-simplicial by Bj\"orner \cite{bjorner:orderings} (see \cite{bjorner.brenti:coxeter} Theorem 3.2.7).  Cambrian lattices, introduced by Nathan Reading, are certain lattice quotients of the weak order \cite{reading:cambrian_lattices}.  Consequently, Theorem \ref{thm_main_quotient}(1) implies that every interval of a Cambrian lattice is either contractible or homotopy equivalent to a sphere.  This property of Cambrian lattices was proved by Reading more generally for any lattice quotient of the chamber poset of a simplicial hyperplane arrangement (\cite{reading:lattice_congruence} Theorem 5.1(iv)).



A Cambrian lattice is an example of a semidistributive lattice.  A lattice $L$ is \emph{meet-semidistributive} if $L$ satisfies
$$x\wedge z=y\wedge z\ \mbox{ implies }\ (x\vee y)\wedge z=x\wedge z\ \mbox{ for }x,y,z\in L.$$
It is known that the M\"obius function of a meet-semidistributive lattice takes values only in $\{1,-1,0\}$ (\cite{mulle:sb} Remark 3.11).  This is a consequence of the following result.

\begin{theorem}\label{thm_main_semi}
If $L$ is a meet-semidistributive lattice, then $L$ is crosscut-simplicial.
\end{theorem}

A lattice is \emph{join-distributive} if it is meet-semidistributive and upper-semimodular, which holds if and only if every atomic interval is Boolean.  Edelman proved that a lattice is join-distributive if and only if it is the dual of some lattice of convex sets of an abstract convex geometry (\cite{edelman:meet} Theorem 3.3).  M{\"u}hle showed that a join-distributive lattice inherits an SB-labelling from this geometry (\cite{mulle:sb} Theorem 3.8).  Distributive lattices were previously shown to inherit an SB-labelling from its subposet of join-irreducibles (\cite{hersh.meszaros:sb} Theorem 5.1).



\begin{figure}
\begin{centering}
\begin{tikzpicture}[scale=.5]
  \draw (-3,-4) -- (3,4)
        (-5,0) -- (5,0)
        (-3,4) -- (3,-4);
  \draw (0,-5/2) node[anchor=north]{$c_0$}
        (0,5/2) node[anchor=south]{$-c_0$};
  \begin{scope}[xshift=10cm]
    \draw (0,-5/2) node[anchor=north]{$c_0$} -- (-4/2,-3/2) -- (-4/2,3/2) -- (0,5/2) node[anchor=south]{$-c_0$} -- (4/2,3/2) -- (4/2,-3/2) -- cycle;
  \end{scope}
\end{tikzpicture}
\caption{\label{fig_ex_chamber}\scriptsize (left) An arrangement of three lines  (right) Poset of chambers}
\end{centering}
\end{figure}
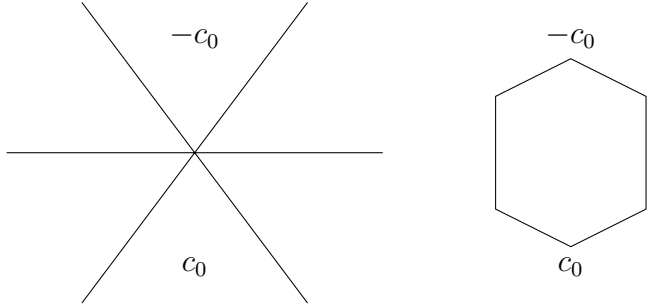

A finite arrangement of linear hyperplanes in $\Rbb^n$ divides the space into a complete fan whose maximal cones are called chambers.  Given an arrangement $\Acal$ with a distinguished chamber $c_0$, the poset $\Pcal(\Acal,c_0)$ is an ordering of the chambers where $c\leq c^{\pr}$ if any hyperplane separating $c_0$ and $c$ also separates $c_0$ and $c^{\pr}$; see figure \ref{fig_ex_chamber}.

The weak order of a finite Coxeter group may be realized as a chamber poset $\Pcal(\Acal,c_0)$ where $\Acal$ is the set of reflecting hyperplanes of the standard reflection arrangement and $c_0$ is any chamber.  The problem of finding other chamber posets that are lattices was considered by Bj\"orner, Edelman, and Ziegler \cite{bjorner.edelman.ziegler:lattice}.  A non-exhaustive list of chamber lattices include posets $\Pcal(\Acal,c_0)$ where
\begin{enumerate}
\item every chamber of $\Acal$ is a simplicial cone (\cite{bjorner.edelman.ziegler:lattice} Theorem 3.4),
\item the intersection lattice of $\Acal$ is supersolvable and $c_0$ is incident to a modular flag of intersection subspaces (\cite{bjorner.edelman.ziegler:lattice} Theorem 4.6),
\item the rank of $\Acal$ is at most 3 and $c_0$ is a simplicial cone (\cite{bjorner.edelman.ziegler:lattice} Theorem 3.2), or
\item the arrangement $\Acal$ is hyperfactored with respect to $c_0$ (\cite{jambu.paris:combinatorics} Theorem 5.1).
\end{enumerate}
The latter three examples may not be crosscut-simplicial, but the first example is crosscut-simplicial.  This follows from Theorem \ref{thm_main_semi} and the fact that $\Pcal(\Acal,c_0)$ is semidistributive if every chamber of $\Acal$ is simplicial (\cite{reading:lattice_hyperplane} Theorem 3).  The semidistributivity property of the weak order was shown previously by Le Conte de Poly-Barbut \cite{poly-barbut:treillis}.


Edelman and Walker proved that every interval of any chamber poset $\Pcal(\Acal,c_0)$ is either contractible or homotopy equivalent to a sphere \cite{edelman.walker:homotopy}.  We characterize those arrangements for which $\Pcal(\Acal,c_0)$ is a crosscut-simplicial lattice in the following Theorem.  We define the bineighborly property in Section \ref{sec_chamber}.


\begin{theorem}\label{thm_main_bineighborly}
Let $\Acal$ be an arrangement with fundamental chamber $c_0$.  The following are equivalent.
\begin{enumerate}
\item $\Pcal(\Acal,c_0)$ is a crosscut-simplicial lattice.
\item $\Pcal(\Acal,c_0)$ is a semidistributive lattice.
\item $\Acal$ is bineighborly.
\end{enumerate}
\end{theorem}

The paper is organized as follows.  In Section \ref{sec_CSL} we recall some background on poset topology and relate SB-labellings with the crosscut-simplicial property.  In Section \ref{sec_semidistributive} we prove Theorem \ref{thm_main_semi} and determine the crosscut complex of a join-semidistributive lattice.  In Section \ref{sec_chamber}, we discuss chamber posets and prove Theorem \ref{thm_main_bineighborly}.  All of the results in this section generalize to the setting of oriented matroids.


We prove part 1 of Theorem \ref{thm_main_quotient} in Section \ref{sec_quotient}.  In Section \ref{sec_doubling} we compute the crosscut complexes of intervals of a lattice doubled at a convex set.  In particular, we prove that every interval of a congruence-normal lattice is either contractible or homotopy equivalent to a sphere.  We also deduce part 2 of Theorem \ref{thm_main_quotient}.

\section{Crosscut-simplicial lattices}\label{sec_CSL}

A \emph{poset} $P$ is a set with a reflexive, anti-symmetric, and transitive binary relation, usually denoted $\leq$.  Unless stated otherwise, all of our posets are assumed to be finite.  The \emph{dual} $(P^*,\leq^*)$ of a poset $(P,\leq)$ is an ordering on the same set where $x\leq^* y$ if and only if $y\leq x$.  Given elements $x,y\in P$ for which $x\leq y$, the \emph{closed interval} $[x,y]$ (\emph{open interval} $(x,y)$) is the set of $z\in P$ such that $x\leq z\leq y$ ($x<z<y$).  An \emph{order filter} (\emph{order ideal}) is a subset $I$ of $P$ such that if $x\in I,\ y\in P$ and $x\geq y$ ($x\leq y$), then $y\in I$.  The \emph{M\"obius function} $\mu$ of a finite poset $P$ is the unique assignment of integers to closed intervals of $P$ for which $\mu([x,y])=1$ if $x=y$ and $\sum_{z\in[x,y]}\mu([x,z])=0$ if $x<y$.  We refer to Chapter 3 of \cite{stanley:enumerative.v1} for applications of M\"obius functions.

A \emph{lattice} is a poset for which every pair of elements $x,y$ has a least upper bound $x\vee y$ and a greatest lower bound $x\wedge y$.  The \emph{join} (\emph{meet}) of a finite subset $A$ of a lattice, denoted $\bigvee A$ ($\bigwedge A$), is the common least upper bound (greatest lower bound) of the elements of $A$.  An atom of an interval $[x,y]$ or $(x,y)$ is any element $z$ covering $x$ with $z\leq y$.  We call an interval $[x,y]$ \emph{atomic} if the join of its atoms is equal to $y$.

Let $(\Delta,A)$ be an abstract simplicial complex on the ground set $A$, and let $F\in\Delta$.  Let $\|\Delta\|$ denote a topological space triangulated by $\Delta$.  The \emph{deletion} $\dl_{\Delta}(F)$ of $F$ is the subcomplex of $\Delta$ of faces disjoint from $F$.  The \emph{star} $\st_{\Delta}(F)$ of $F$ is the subcomplex of faces $F^{\pr}$ such that $F\cup F^{\pr}\in\Delta$.  The \emph{link} $\lk_{\Delta}(F)$ of $F$ is the subcomplex of $\st_{\Delta}(F)$ of faces disjoint from $F$.

The \emph{join} $\Delta*\Delta^{\pr}$ of two complexes $(\Delta,A),(\Delta^{\pr},A^{\pr})$ is the simplicial complex on $A\sqcup A^{\pr}$ with faces $F\sqcup F^{\pr}$ where $F\in\Delta,\ F^{\pr}\in\Delta^{\pr}$.  The join of abstract simplicial complexes realizes the topological join $\|\Delta*\Delta^{\pr}\|\cong\|\Delta\|*\|\Delta^{\pr}\|$.  The \emph{cone} $\{v\}*\Delta$ is the join of $\Delta$ with a one-element complex.  The \emph{suspension} $\{v,v^{\pr}\}*\Delta$ is the join of $\Delta$ with a discrete two-element complex.

The \emph{order complex} of a poset $P$ is the simplicial complex of chains $x_0<\cdots<x_d$ of elements of $P$.  The order complex is self-dual, even if the poset is not.  If $P$ is the set of faces of a simplicial complex $X$ ordered by inclusion, then the order complex of $P$ is homeomorphic to $X$.  Thus, we define the topology of a poset to be that of its order complex.  The link of a face $x_0<\cdots<x_d$ is isomorphic to the join of the order complexes of $P_{<x_0},(x_0,x_1),\ldots,(x_{d-1},x_d),P_{>x_d}$.  Hence, the local topology of $P$ is completely determined by the topology of intervals and principal order ideals and filters of $P$.

We write $\hat{1}$ ($\hat{0}$) for the largest (smallest) element of a poset $P$, if it exists.  A poset is \emph{bounded} if it contains a largest element and a smallest element.  The \emph{proper part} $\ov{P}$ of a bounded poset $P$ is the same poset with those bounds removed.

If $P$ is a bounded poset, the reduced Euler characteristic of $\ov{P}$ is equal to $\mu(\hat{0},\hat{1})$.  Hence, $\mu(\hat{0},\hat{1})$ is a homotopy invariant of $\ov{P}$.  The full M\"obius function is then determined by the local topology of $P$.  Many methods for computing homotopy invariants of posets are given in Section 10 of Bj\"orner's survey \cite{bjorner:topological}.  We only require the Crosscut Theorem.

A crosscut $C$ of a poset $P$ is a subset of pairwise incomparable elements satisfying the following two conditions.
\begin{itemize}
\item For every chain $x_0<\cdots<x_d$ of $P$, there exists an element of $C$ comparable to every $x_i$.
\item If $B\subseteq C$, then $B$ has either at most one common minimal upper bound or one common maximal lower bound.
\end{itemize}

If $C$ is a crosscut of $P$, then the \emph{crosscut complex} $\Gamma(P,C)$ is the simplicial complex on $C$ containing subsets $B\subseteq C$ for which either $\bigvee B$ or $\bigwedge B$ exists.

\begin{theorem}[Crosscut Theorem]
If $C$ is a crosscut of $P$, then $P$ is homotopy equivalent to its crosscut complex $\Gamma(P,C)$.
\end{theorem}

The set of atoms of a lattice is always a crosscut.  Unless specified otherwise, we will refer to this crosscut when discussing the crosscut complex of a lattice.

We say a lattice is \emph{crosscut-simplicial} if for every interval $(x,y)$, the crosscut complex on the atoms $A$ of $(x,y)$ is either a $(|A|-1)$-simplex or the boundary of a $(|A|-1)$-simplex.  This notion is not self-dual as the crosscut complex on the coatoms of $(x,y)$ may not be isomorphic to the crosscut complex on the atoms.  It is immediate from the definitions that a lattice is crosscut-simplicial if and only if for every interval $[x,y]$, the top element $y$ is not the join of any proper subset of atoms of $[x,y]$.

An \emph{SB-labelling} of a lattice is a labelling $\lambda$ of the covering relations such that
\begin{list}{}{}
\item[(SB1)] if $y$ and $z$ are distinct elements covering some element $x$, then $\lambda(x\lessdot y)$ is distinct from $\lambda(x\lessdot z)$; and
\item[(SB2)] if $B$ is a subset of atoms of $(x,\hat{1})$, then every saturated chain from $x$ to $\bigvee B$ contains only labels in the set $\{\lambda(x\lessdot y):\ y\in B\}$, and each of those labels occurs at least once.
\end{list}

The topological significance of an SB-labelling is encapsulated in the following theorem.

\begin{theorem}[Hersh-M\'esz\'aros \cite{hersh.meszaros:sb} Theorem 3.7]
If $L$ admits a SB-labelling, then $L$ is crosscut-simplicial.
\end{theorem}

\begin{figure}
\begin{centering}
\begin{tikzpicture}
\draw (0,0) -- (3,2.5)
    (0,1) -- (-1,2)
    (0,2) -- (-1,3)
    (1,2) -- (0,3)
    (1,3) -- (0,4);
\draw[densely dashed] (0,0) -- (0,1) -- (0,2)
    (-1,2) -- (-1,3)
    (1,2) -- (1,3)
    (0,3) -- (0,4) -- (3,2.5);
\draw[->,decorate,decoration={snake,amplitude=.4mm,segment length=4mm,post length=1mm}] (0,2) -- (1,3);
\draw[->,decorate,decoration={snake,amplitude=.4mm,segment length=4mm,post length=1mm}] (0,3) -- (-1,2);
\draw[->,decorate,decoration={snake,amplitude=.4mm,segment length=4mm,post length=1mm}] (0,1) -- (1,2);
\draw[->,decorate,decoration={snake,amplitude=.4mm,segment length=4mm,post length=1mm}] (0,4) -- (-1,3);
\end{tikzpicture}
\caption{\label{fig_sb_labelling}\scriptsize A lattice with an edge-labelling satisfying (SB1) and (SB2$^{\pr}$).}
\end{centering}
\end{figure}
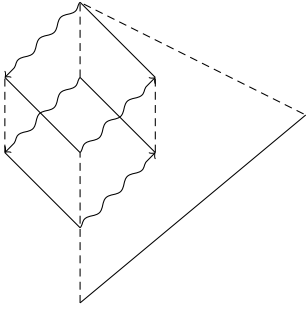

The converse need not hold, as shown in Figure \ref{fig_sb_labelling}.  However, this example does have an SB-labelling if one relaxes condition (2) by
\begin{list}{}{}
\item[(SB2$^{\pr}$)] if $B$ is a subset of atoms $A$ of $(x,\hat{1})$, then every saturated chain from $x$ to $\bigvee B$ contains each label in the set $\{\lambda(x\lessdot y):\ y\in B\}$ at least once, and it contains no labels from the set $\{\lambda(x\lessdot z):\ z\in A-B\}$.
\end{list}

Lattices with an edge labelling satisfying (SB1) and (SB2$^{\pr}$) are crosscut-simplicial.  We do not know of a crosscut-simplicial lattice that does not admit this more general version of SB-labelling.

\section{Semidistributive Lattices}\label{sec_semidistributive}

A lattice $L$ is \emph{meet-semidistributive} if for $x,y,z\in L$, the equality $x\wedge z=y\wedge z$ implies $(x\vee y)\wedge z=x\wedge z$.  It is \emph{join-semidistributive} if it satisfies a dual condition.  A lattice is \emph{semidistributive} if it is both meet-semidistributive and join-semidistributive.

\begin{theorem}\label{thm_meet_semi}
If $L$ is a meet-semidistributive lattice, then $L$ is crosscut-simplicial.
\end{theorem}

\begin{poof}
Let $A$ be the set of atoms of $L$.  Since meet-semidistributivity is inherited by intervals, it suffices to prove $\bigvee B<\bigvee A$ whenever $B$ is a proper subset of $A$.

We proceed by induction on $|A|$.  Let $B$ be a minimal subset of $A$ such that $\bigvee B=\bigvee A$, and let $x\in B$.  If $A-x$ contains an element $z$ such that $z$ and $\bigvee(B-x)$ are incomparable, then $x\wedge z=\hat{0}=\bigvee(B-x)\wedge z$ holds.  But this implies $z=(x\vee\bigvee(B-x))\wedge z=\hat{0}$, a contradiction.  Hence, $A-x$ is the set of atoms of the meet-distributive lattice $[\hat{0},\bigvee(B-x)]$.  By induction, we have $A-x=B-x$, as desired.
\end{poof}

\begin{corollary}\label{cor_join_semi}
Every interval of a meet-semidistributive lattice or join-semidistributive lattice is either contractible or homotopy equivalent to a sphere.
\end{corollary}

Figure \ref{fig_sb_labelling} shows an example of a crosscut-simplicial lattice that is not meet-semidistributive.  A join-semidistributive lattice may not be crosscut-simplicial, but its crosscut complex still admits a simple discription.

\begin{proposition}\label{thm_join_semi}
If $L$ is a join-semidistributive lattice with atom set $A$, then its crosscut complex is either a $(|A|-1)$-simplex or a pure $(|A|-2)$-subcomplex of the $(|A|-1)$-simplex.
\end{proposition}

\begin{poof}
If $\bigvee A<\hat{1}$, then the crosscut complex of $L$ is a $(|A|-1)$-simplex.  Hence, we may assume $\bigvee A=\hat{1}$.  We prove that the maximal faces of the crosscut complex are all of dimension $|A|-2$.

Let $B$ be a maximal subset of $A$ such that $\bigvee B<\hat{1}$.  Suppose $A-B$ has two distinct elements $x,y$.  By the maximality of $B$, one has $x\vee(\bigvee B)=\hat{1}=y\vee(\bigvee B)$.  But this implies $(x\wedge y)\vee(\bigvee B)=\hat{1}$, which is impossible since $(x\wedge y)\vee(\bigvee B)=\bigvee B<\hat{1}$.
\end{poof}

\section{Poset of chambers of a hyperplane arrangement}\label{sec_chamber}

A \emph{real, central hyperplane arrangement} $\Acal$ is a finite set of hyperplanes in $\Rbb^n$ whose common intersection contains the origin.  The arrangement determines a complete fan of cones, the \emph{faces} of the arrangement, whose maximal cones are called \emph{chambers}.  A \emph{wall} of a chamber $c$ is any hyperplane in $\Acal$ incident to $c$.  We let $\Tcal(\Acal)$ denote the set of chambers of $\Acal$ and $\Wcal(c)$ the walls of a chamber $c$.  Given two chambers $c,c^{\pr}\in\Tcal(\Acal)$, the \emph{separation set} $S(c,c^{\pr})$ is the set of hyperplanes in $\Acal$ separating $c$ and $c^{\pr}$.  Given a chamber $c_0\in\Tcal(\Acal)$, the \emph{poset of chambers} $\Pcal(\Acal,c_0)$ is an ordering on $\Tcal(\Acal)$ where $c\leq c^{\pr}$ if $S(c_0,c)\subseteq S(c_0,c^{\pr})$.  The distinguished chamber $c_0$ is called the \emph{fundamental chamber}.  The \emph{intersection lattice} $L(\Acal)$ is the set of intersection subspaces of $\Acal$ ordered by reverse inclusion.  If $X\in L(\Acal)$, the \emph{localization} $\Acal_X$ is the subarrangement of hyperplanes containing $X$.  If $\Acal^{\pr}\subseteq\Acal$, the \emph{restriction} $c|_{\Acal^{\pr}}$ of a chamber $c$ is the unique chamber of $\Acal^{\pr}$ containing $c$.  A chamber is \emph{simplicial} if its face poset is Boolean.

The set of chambers $\Tcal(\Acal)$ determines the oriented matroid structure of $\Acal$ (see e.g. \cite{bjorner.lasVergnas.ea:oriented} Proposition 3.8.2), while the intersection lattice $L(\Acal)$ determines its underlying matroid.  All of the results in this section hold in the more general context of oriented matroids; see Section 2.1 of \cite{bjorner.lasVergnas.ea:oriented} for the translation.  As we will not need this level of generality, we stick with the language of hyperplane arrangements.

We collect several key properties of chamber posets in the following proposition.

\begin{proposition}[\cite{edelman:partial}]\label{prop_chamber_poset}
Let $\Acal$ be an arrangement with a fundamental chamber $c_0$.
\begin{enumerate}
\item\label{prop_chamber_poset_involution} $\Pcal(\Acal,c_0)$ admits a free, order-reversing involution $c\mapsto -c$ such that $\Acal$ is the disjoint union of $S(c_0,c)$ and $S(c_0,-c)$ for any chamber $c$.
\item\label{prop_chamber_poset_graded} $\Pcal(\Acal,c_0)$ is a bounded, graded poset with rank function $c\mapsto|S(c_0,c)|$.
\item\label{prop_chamber_poset_incidence} For $c\in\Tcal(\Acal),\ X\in L(\Acal)$, if $c$ is incident to $X$, then there exists a chamber $c^{\pr}$ such that $S(c,c^{\pr})=\Acal_X$.
\item\label{prop_chamber_poset_opposite} For $c,c^{\pr}\in\Tcal(\Acal)$, if $\Wcal(c)\subseteq S(c,c^{\pr})$ then $c^{\pr}=-c$.
\item\label{prop_chamber_poset_restriction} For $\Acal^{\pr}\subseteq\Acal$, the map $c\mapsto c|_{\Acal^{\pr}}$ is an order-preserving map of posets $\Pcal(\Acal,c_0)\ra\Pcal(\Acal^{\pr},(c_0)|_{\Acal^{\pr}})$.
\end{enumerate}
\end{proposition}

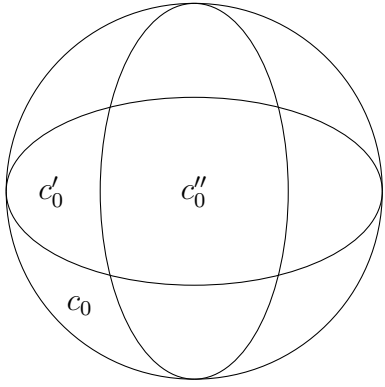
\begin{figure}
\begin{centering}
\begin{tikzpicture}[scale=.5]
\begin{scope}
  \draw (0,0) circle (5cm);
  \draw (5cm,0) arc (0:180:5cm and 2.5cm);
  \draw[rotate=90] (5cm,0) arc (0:180:5cm and 2.5cm);
  \draw[rotate=180] (5cm,0) arc (0:180:5cm and 2.5cm);
  \draw[rotate=270] (5cm,0) arc (0:180:5cm and 2.5cm);
  \draw (225:4.3cm) node{$c_0$};
  \draw (180:3.8cm) node{$c_0^{\pr}$};
  \draw (0,0) node{$c_0^{\pr\pr}$};
\end{scope}
\end{tikzpicture}
\caption{\label{fig_chamber_poset}\scriptsize $\Pcal(\Acal,c_0)$ is semidistributive and $(\Acal,c_0)$ is bineighborly. $\Pcal(\Acal,c_0^{\pr})$ is a non-semidistributive lattice.  $(\Acal,c_0^{\pr})$ is not bineighborly since $c_0^{\prpr}$ has two upper walls that do not intersect at the boundary of $c_0^{\prpr}$.  $\Pcal(\Acal,c_0^{\pr\pr})$ is not a lattice.}
\end{centering}
\end{figure}

The \emph{upper walls} $U(c)$ of a chamber $c$ is the set of hyperplanes $H\in\Acal$ incident to $c$ such that $H\notin S(c_0,c)$.  We say an arrangement $\Acal$ with fundamental chamber $c_0$ is \emph{bineighborly} if for $c\in\Tcal(\Acal),\ H,H^{\pr}\in U(c)$, the chamber $c$ is incident to $H\cap H^{\pr}$.  Since chamber posets are self-dual by Proposition \ref{prop_chamber_poset}(\ref{prop_chamber_poset_involution}), if $\Acal$ is bineighborly, a similar condition holds for the lower walls of any chamber.  The bineighborly property typically depends on the choice of fundamental chamber; see Figure \ref{fig_chamber_poset}.

\begin{theorem}\label{thm_bineighborly}
Let $\Acal$ be a real, central hyperplane arrangement with fundamental chamber $c_0$.  The following are equivalent.
\begin{enumerate}
\item $\Pcal(\Acal,c_0)$ is crosscut-simplicial.
\item $\Acal$ is bineighborly.
\item $\Pcal(\Acal,c_0)$ is a semidistributive lattice.
\end{enumerate}
\end{theorem}

We use the following results in our proof.

\begin{lemma}[\cite{bjorner.edelman.ziegler:lattice} Lemma 2.1]\label{lem_local_lattice}
Let $L$ be a finite bounded poset.  If the join $x\vee y$ exists whenever $x,y$ cover a common element, then $L$ is a lattice.
\end{lemma}

\begin{lemma}\label{lem_chamber_restriction}
For $c\in\Tcal(\Acal)$, the face poset of $c$ is isomorphic to the face poset of $c|_{\Wcal(c)}$.
\end{lemma}

\begin{poof}(of Theorem \ref{thm_bineighborly})
We show $(1)\Ra(2)$ and $(2)\Ra(3)$.  The implication $(3)\Ra(1)$ is a special case of Theorem \ref{thm_meet_semi}.  


$(1)\Ra(2)$: Suppose $\Pcal(\Acal,c_0)$ is crosscut-simplicial.  Let $H$ and $H^{\pr}$ be upper walls of a chamber $c$.  Let $d,d^{\pr}\in\Tcal(\Acal)$ with $S(c,d)=\{H\},\ S(c,d^{\pr})=\{H^{\pr}\}$.  Since $\Pcal(\Acal,c_0)$ is crosscut-simplicial, $d$ and $d^{\pr}$ are the only atoms of the interval $[c,d\vee d^{\pr}]$.  Thus, the set $S(c,d\vee d^{\pr})$ contains no walls of $c$ besides $H$ and $H^{\pr}$.  Let $\alpha$ be a generic point in the intersection of the cones $c|_{\Wcal(c)}$ and $d|_{\Wcal(c)}$, and let $\beta$ be a generic point in the intersection of $d^{\pr}|_{\Wcal(c)}$ and $(d\vee d^{\pr})|_{\Wcal(c)}$.  Both $\alpha$ and $\beta$ are points in $H$ separated by $H^{\pr}$ and no other wall of $c$.  The line segment between $\alpha$ and $\beta$ intersects $H^{\pr}$, so there is a point in $H\cap H^{\pr}$ on the same side as $c$ of any $H^{\prpr}\in\Wcal(c)-\{H,H^{\pr}\}$.  Hence, $c|_{\Wcal(c)}$ is incident to $H\cap H^{\pr}$.  Since the face posets of $c$ and $c|_{\Wcal(c)}$ are the same, the chamber $c$ is incident to $H\cap H^{\pr}$.



$(2)\Ra(3)$: Assume $\Acal$ is bineighborly.  Let $c\in\Tcal(\Acal),\ H,H^{\pr}\in U(c)$, and let $a,b$ be the chambers with $S(c,a)=\{H\},\ S(c,b)=\{H^{\pr}\}$.  By the bineighborly assumption, $c$ is incident to $H\cap H^{\pr}$.  Hence, by Proposition \ref{prop_chamber_poset}(\ref{prop_chamber_poset_incidence}), there exists a chamber $c^{\pr}$ such that $S(c,c^{\pr})=\Acal_{H\cap H^{\pr}}$.  If $d$ is some chamber such that $H,H^{\pr}\in S(c,d)$, then $\Acal_{H\cap H^{\pr}}\subseteq S(c,d)$.  Therefore, $c^{\pr}$ is the join of $a$ and $b$.  By Lemma \ref{lem_local_lattice}, this implies $\Pcal(\Acal,c_0)$ is a lattice.

It remains to prove the following claim.

Claim: For $a\leq b$ and $x,y\in[a,b],\ z\in\Pcal(\Acal,c_0)$, if $x\wedge z=y\wedge z$, then $(x\vee y)\wedge z=x\wedge z$.

If $a=b$ or $a\lessdot b$, the claim is trivial.  Let $a<b$ and suppose the claim holds for all proper subintervals of $[a,b]$.  Let $x,y,z$ be chambers such that $x\wedge z=y\wedge z,\ x,y\in[a,b]$.  We may assume $a=x\wedge y$ and $x\vee y=b$ by the inductive hypothesis.  We have
$$x\wedge z=x\wedge(x\wedge z)=x\wedge(y\wedge z)=a\wedge z.$$

Let $u$ be a coatom of $[x,b]$.  Since $x\leq u<b=x\vee y$ holds, $u$ is not an upper bound for $y$.  We have $a=x\wedge y\leq u\wedge y$.  Since $x\vee(u\wedge y)\leq u<b$ and $u\wedge y\in[a,y]$, the inductive hypothesis implies $(x\vee(u\wedge y))\wedge z=x\wedge z$.  If $a<u\wedge y$ then the inductive hypothesis implies $((x\vee(u\wedge y))\vee y)\wedge z=x\wedge z$.  This simplifies to $(x\vee y)\wedge z=x\wedge z$, as desired.  Thus, we may assume $a=u\wedge y$ holds for every coatom $u$ of $[x,b]$.

Let $d$ be an atom of $[a,y]$.  If $u$ is a coatom of $[x,b]$, then $d$ and $u$ are incomparable, so $S(a,d)=S(u,b)=\{H\}$ for some hyperplane $H$.  In particular, $H$ is the unique hyperplane in $S(a,y)\cap U(a)$ and in $S(x,b)\cap U(b)$.

Let $c$ be an atom of $[a,x]$ and let $H^{\pr}$ be the hyperplane separating $c$ and $a$.  By assumption, $a$ is incident to $H\cap H^{\pr}$, so by Proposition \ref{prop_chamber_poset}(\ref{prop_chamber_poset_incidence}) the join of $c$ and $d$ is the chamber satisfying $S(a,c\vee d)=\Acal_{H\cap H^{\pr}}$.

Since $c\in[a,x],\ d\in[a,y]$, we have 
$$c\wedge z=a\wedge z=d\wedge z.$$

Assume $c\vee d=b$.  Then $S(a\wedge z,b\wedge z)\subseteq S(a,b)=\Acal_{H\cap H^{\pr}}$ holds.  Neither $H$ nor $H^{\pr}$ is in $S(a\wedge z,b\wedge z)$ since $c\wedge z=a\wedge z=d\wedge z$.  Hence, $a\wedge z=b\wedge z$.

Now assume $c\vee d<b$.  Since $x$ and $c\vee d$ are both in $[c,b]$, the equality $x\wedge z=(x\vee(c\vee d))\wedge z$ holds by induction.  Similarly, $y\wedge z=(y\vee c\vee d)\wedge z$.  Finally, $x\vee(c\vee d)$ and $y\vee(c\vee d)$ are both elements of $[c\vee d,b]$, so
$$x\wedge z=((x\vee(c\vee d))\vee(y\vee(c\vee d)))\wedge z=(x\vee y)\wedge z.$$
\end{poof}

\section{Lattice Quotients}\label{sec_quotient}

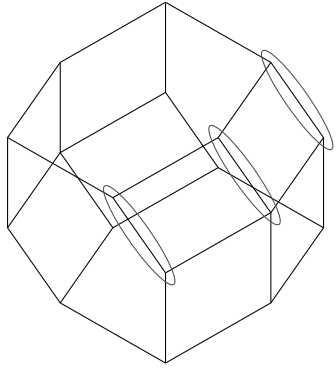
\begin{figure}
\begin{centering}
\begin{tikzpicture}
\coordinate (A) at (-1.4,.8);
\coordinate (B) at (-.7,1);
\coordinate (Bp) at (-.35,.5);
\coordinate (C) at (0,1.2);
\coordinate (D) at (.7,1);
\coordinate (E) at (1.4,.8);

\draw (0,0) -- (A) -- ($(A)+(B)$) -- ($(A)+(B)+(C)$) -- ($(A)+(B)+(C)+(D)$) -- ($(A)+(B)+(C)+(D)+(E)$) -- ($(B)+(C)+(D)+(E)$) -- ($(C)+(D)+(E)$) -- ($(D)+(E)$) -- (E) -- cycle;
\draw (0,0) -- (C) -- ($(B)+(C)$) -- ($(A)+(B)+(C)$);
\draw ($(B)+(C)$) -- ($(B)+(C)+(E)$) -- ($(B)+(C)+(D)+(E)$);
\draw (C) -- ($(C)+(E)$) -- ($(B)+(C)+(E)$);
\draw (E) -- ($(C)+(E)$) -- ($(C)+(D)+(E)$);
\draw ($(A)+(B)+(C)+(D)+(E)$) -- ($(A)+(B)+(D)+(E)$) -- ($(A)+(D)+(E)$) -- ($(D)+(E)$);
\draw ($(A)+(D)+(E)$) -- ($(A)+(D)$) -- (A);
\draw ($(A)+(B)+(D)+(E)$) -- ($(A)+(B)+(D)$) -- ($(A)+(D)$);
\draw ($(A)+(B)+(C)+(D)$) -- ($(A)+(B)+(D)$) -- ($(A)+(B)$);
\draw[black!70] ($(C)+(Bp)$) ellipse[rotate=35, x radius=.15cm, y radius=.8cm];
\draw[black!70] ($(C)+(E)+(Bp)$) ellipse[rotate=35, x radius=.15cm, y radius=.8cm];
\draw[black!70] ($(C)+(D)+(E)+(Bp)$) ellipse[rotate=35, x radius=.15cm, y radius=.8cm];
\end{tikzpicture}
\caption{\label{fig_quotient}\scriptsize A lattice congruence on $\Pcal(\Acal,c_0)$ of Figure \ref{fig_chamber_poset}.}
\end{centering}
\end{figure}

A \emph{lattice congruence} $\Theta$ is an equivalence relation on a lattice $L$ respecting the meet and join operations; see Figure \ref{fig_quotient}.  More precisely, for $x,y,z\in L$ if $x\equiv y(\Theta)$ holds then so does $x\vee z\equiv y\vee z(\Theta)$ and $x\wedge z\equiv y\wedge z(\Theta)$.  If $L$ is finite, then the equivalence classes are necessarily closed intervals.  As a result, there are order-preserving maps $\pi^{\uparrow},\pi^{\downarrow}$ on $L$ where $\pi^{\uparrow}(x)$ ($\pi^{\downarrow}(x)$) is the largest (smallest) element in the $\Theta$-equivalence class $[x]$.

Given a lattice congruence $\Theta$, the \emph{quotient lattice} $L/\Theta$ is the collection of equivalence classes with $[x]\leq[y]$ if there exists some elements $x^{\pr}\in[x]$ and $y^{\pr}\in[y]$ such that $x^{\pr}\leq y^{\pr}$.  The images of $\pi^{\uparrow}$ and $\pi^{\downarrow}$ are both isomorphic to $L/\Theta$.

We summarize the key properties of quotient lattices described in Section 2 of \cite{reading:lattice_congruence}.

\begin{lemma}[Reading \cite{reading:lattice_congruence} Lemma 2.1]\label{lem_interval_restriction}
Let $\Theta$ be a lattice congruence of $L$ and $[x,y]$ an interval of $L$.  The restriction of $\Theta$ to $[x,y]$ is a lattice congruence of $[x,y]$.  Moreover, the interval $[[x],[y]]$ of $L/\Theta$ is isomorphic to $[x,y]/\Theta$.
\end{lemma}

\begin{lemma}[Reading \cite{reading:lattice_congruence} Proposition 2.2]\label{lem_atoms}
The atoms of $L/\Theta$ are in bijection with the set of elements covering $\pi^{\uparrow}(\hat{0})$ via the map $a\mapsto[a]$.
\end{lemma}

\begin{poof}
Let $x=\pi^{\uparrow}(\hat{0})$.  We verify that the above map $a\mapsto[a]$ is a well-defined, bijective map from covers of $x$ to covers of $[x]$.

Let $a\in L$ be a cover of $x$.  If $[b]<[a]$, then $[b\vee x]=[b]$, which implies $\pi^{\downarrow}(b)\vee x<a$.  Since $a$ covers $x$, this forces $[b]=[\hat{0}]$.

Assume $a,b$ cover $x$ such that $[a]=[b]$.  Since $[a\wedge b]=[a],\ [a]\neq[x]$, we have $x<a\wedge b\leq a$.  But $a$ covers $x$, so $a\leq b$.  Similarly, $b\leq a$.

Assume $[a]$ covers $[x]$ and let $a^{\pr}$ be the smallest element in the class $[a]$ larger than $x$.  If $x<b<a^{\pr}$ for some $b\in L$, then $[x]<[b]<[a]$, an impossibility.  Hence, $a^{\pr}$ covers $x$.
\end{poof}

The main theorem of this section is essentially a restatement of Corollary 2.4 of \cite{reading:lattice_congruence}.

\begin{theorem}\label{thm_crosscut_quotient}
The crosscut complex of any interval of $L/\Theta$ is isomorphic to the crosscut complex of some interval of $L$.
\end{theorem}

\begin{poof}
Let $([x],[y])$ be an interval of $L/\Theta$.  Let $A$ be the set of atoms of $(\pi^{\uparrow}(x),\pi^{\uparrow}(y))$.  Let $y^{\pr}$ be the smallest element $\Theta$-equivalent to $y$ such that $\bigvee A\leq y^{\pr}$.  We claim that the crosscut complex of $(\pi^{\uparrow}(x),y^{\pr})$ is isomorphic to that of $([x],[y])$.

If $\bigvee A<y^{\pr}$, then $[\bigvee A]<[y]$, and both complexes are isomorphic to a $(|A|-1)$-simplex.  Thus, we may assume $\bigvee A=y^{\pr}$.

By Lemma \ref{lem_atoms}, the map $a\mapsto[a]$ is a bijection on the sets of atoms of $(\pi^{\uparrow}(x),y^{\pr})$ and $([x],[y])$.  Let $B\subseteq A$.  If $\bigvee B=y^{\pr}$, then $\bigvee_{b\in B}[b]=[\bigvee B]=[y]$.

Conversely, suppose $\bigvee_{b\in B}[b]=[y]$ and assume $\bigvee B<y^{\pr}$.  Then there exists $a\in A-B$ such that $\bigvee B<a\vee(\bigvee B)$ since $\bigvee B<\bigvee A$.  Since $a$ covers $\pi^{\uparrow}(x)$, this forces $a\wedge(\bigvee B)=\pi^{\uparrow}(x)$.  But, $[x]<[a]<[\bigvee B]$ so $[a]\wedge[\bigvee B]\neq[x]$, a contradiction.
\end{poof}

\begin{corollary}
Let $L$ be a crosscut-simplicial lattice.  If $\Theta$ is a lattice congruence of $L$, then the quotient $L/\Theta$ is crosscut-simplicial.
\end{corollary}


\section{Doubling}\label{sec_doubling}

A subset $C$ of a poset $P$ is \emph{order-convex} if for $x,y\in C,\ x\leq y,$ the interval $[x,y]$ is contained in $C$.  If $C$ is a subset of $P$, we let $P_{\geq C}$ be the subposet of elements $x\in P$ for which there exists $c\in C$ with $x\geq c$.  The \emph{doubling} $P[C]$ of $P$ at an order-convex subset $C$ is the induced subposet of $P\times\{0,1\}$ on the set
$$P[C]=((P-P_{\geq C})\cup C)\times\{0\}\ \cup\ P_{\geq C}\times\{1\}.$$
It is straight-forward to check that $L[C]$ is a lattice if $L$ is a lattice.

\begin{lemma}\label{lem_double_lattice}
Let $C$ be an order-convex subset of a lattice $L$.  If $(x,\epsilon)$ and $(y,\epsilon^{\pr})$ are elements of $L[C]$, then their join is
$$(x,\epsilon)\vee(y,\epsilon^{\pr})=\begin{cases}(x\vee y,\max\{\epsilon,\epsilon^{\pr}\})\ &\mbox{if }x\vee y\in (P-P_{\geq C})\cup C\\(x\vee y,1)\ &\mbox{otherwise.}\end{cases}$$
\end{lemma}

The projection $\pi:P[C]\ra P$ defined by $\pi(x,\epsilon)=x$ is order-preserving.  If $P$ is a lattice, the map $\pi$ is a lattice quotient map (see Section \ref{sec_quotient}).  The crosscut complexes of a doubled lattice are related to those of the original lattice as in the following proposition.  We let $\Delta(A)$ denote the simplicial complex of all subsets of $A$.  If $\Gamma$ is a simplicial complex and $B$ a subset of the ground set, we let $\Gamma|_B$ denote the induced subcomplex on $B$.

\begin{proposition}\label{prop_doubling_crosscut}
Let $C$ be an order-convex subset of a lattice $L$.  Let $I$ be an open interval of $L[C]$ and let $A$ be the set of atoms of $\pi(I)$.  At least one of the following holds.
\begin{enumerate}
\item $\Gamma(I)\cong\Gamma(\pi(I))$
\item $\Gamma(I)\cong\Delta(A)$
\item $\Gamma(I)\cong\{v\}*\Gamma(\pi(I))|_{A\cap C}$
\item $\Gamma(I)\cong\Delta(A)\cup(\{v\}*\Gamma(\pi(I)))$
\end{enumerate}
\end{proposition}

\begin{poof}(of Proposition \ref{prop_doubling_crosscut})
We divide the possible intervals of $L[C]$ into four cases.  Let $I=((x,\epsilon),(y,\epsilon^{\pr}))$ be an interval.  Then exactly one of the following holds:
\begin{enumerate}
\item $\epsilon=\epsilon^{\pr}$ or $x,y\notin C$,
\item $\epsilon<\epsilon^{\pr},\ x\notin C$, and $y\in C$,
\item $\epsilon<\epsilon^{\pr},\ x\in C$, and $y\notin C$, or
\item $\epsilon<\epsilon^{\pr}, x\in C$, and $y\in C$.
\end{enumerate}
We verify that these line up with the four cases for $\Gamma(I)$ listed above.

(1) If $\epsilon=\epsilon^{\pr}$, then the open interval $((x,\epsilon),(y,\epsilon^{\pr}))$ is isomorphic to $(x,y)$.  If $\epsilon<\epsilon^{\pr}$ and $x,y$ are both not in $C$, then by Lemma \ref{lem_double_lattice} the join of some atoms $B$ in $((x,0),(y,1))$ is equal to $(y,1)$ if and only if the join of $\{\pi(b):\ b\in B\}$ equals $y$.  In both cases, the crosscut complexes of $I$ and $\pi(I)$ are isomorphic.

(2) If $\epsilon<\epsilon^{\pr},\ x\notin C$, and $y\in C$, then the join of all of the atoms of $I$ is bounded above by $(y,0)$.  Hence, $\Gamma(I)$ is a simplex.

(3) Suppose $\epsilon<\epsilon^{\pr},\ x\in C$, and $y\notin C$.  Then
$$\{(x,1)\}\cup\{(a,0):\ a\in A\cap C\}$$
is the set of atoms of $I$.  If $B$ is a set of atoms of $I$ whose join is equal to $(z,\epsilon^{\prpr})$, then $\bigvee(B\cup\{(x,1)\})$ equals $(z,1)$.  If $(z,\epsilon^{\prpr})<(y,1)$ then $(z,1)<(y,1)$ since $y\notin C$.  If $(z,\epsilon^{\prpr})=(y,1)$, then $\bigvee_{b\in B}\pi(b)=y$.  Hence, $\Gamma(I)$ is the cone $\{(x,1)\}*\Gamma(\pi(I))|_{A\cap C}$.

(4) If $\epsilon<\epsilon^{\pr}, x\in C$, and $y\in C$, then $I$ is isomorphic to $\ov{[x,y]\times 2}$.  Let $\phi:A\ra I$ be the inclusion $a\mapsto(a,0)$.  The atom set of $I$ is $\{(x,1)\}\cup\phi(A)$.  Since $(a,0)\leq (y,0)$ for $a\in A$, the deletion $\Gamma(I)-\{(x,1)\}$ is equal to $\Delta(\phi(A))$.  If $A^{\pr}\subseteq A$ then $(x,1)\vee\bigvee\phi(A)$ equals $(\bigvee A,1)$, so the link of $(x,1)$ is equal to $\phi(\Gamma(\pi(I)))$.
\end{poof}

The set $\Con(L)$ of congruences of a lattice $L$ ordered by inclusion forms a distributive lattice \cite{funayama.nakayama:distributivity}.  In particular, for any cover relation $x\lessdot y$, there is a minimal congruence $\Cg(x,y)$ of $L$ in which $x$ and $y$ are equivalent.  If $x\in L$ is meet-irreducible (join-irreducible) there is a unique element $x^*$ covering $x$ ($x_*$ covered by $x$).  A lattice $L$ is \emph{congruence-normal} if for any meet-irreducible $x$ and join-irreducible $y$ in $L$, the equality $\Cg(x,x^*)=\Cg(y_*,y)$ implies $x\ngeq y$.  Some interesting examples of congruence-normal lattices are distributive lattices, the weak order of a finite Coxeter group (\cite{caspard.poly-barbut.morvan:cayley} Theorem 6 or \cite{reading:lattice_hyperplane} Theorem 27), and chamber posets of supersolvable arrangements (\cite{reading:lattice_hyperplane} Theorem 1).  As Cambrian lattices are quotients of the weak order, they also inherit the congruence-normal property.

Day proved that a lattice is congruence-normal if and only if it can be obtained from the one-element poset by a sequence of doublings (\cite{day:congruence} Section 3).  Using Proposition \ref{prop_doubling_crosscut}, we deduce the corollary.

\begin{corollary}\label{cor_congruence_normal}
Let $L$ be a lattice with an order convex subset $C$.  If every interval of $L$ is either contractible or homotopy equivalent to a sphere, then the same holds for $L[C]$.  In particular, if $L$ is congruence-normal, then every interval of $L$ is either contractible or homotopy equivalent to a sphere.
\end{corollary}

\begin{poof}
Let $I$ be a closed interval of $L$.  It suffices to show that the four constructions for $\Gamma(\ov{I})$ in Proposition \ref{prop_doubling_crosscut} preserve the property of being contractible or homotopy equivalent to a sphere.  The first two cases are trivial.  The third case is a cone, so it is contractible.  In the fourth case, $I$ is isomorphic to $\pi(I)\times 2$, so $\ov{I}$ is homeomorphic to the suspension of $\ov{\pi(I)}$ by Theorem 5.1(d) of \cite{walker:homeomorphisms}.
\end{poof}

\begin{figure}
\begin{centering}
\begin{tikzpicture}

\begin{scope}
\filldraw (0,0) circle(.5mm);
\draw[black!70] (0,0) circle(2mm);
\end{scope}

\begin{scope}[xshift=.8cm,yshift=-.75cm]
\draw (0,0) -- (0,1.5);
\draw[black!70] (0,.75) ellipse[x radius=.2cm, y radius=.9cm];
\end{scope}

\begin{scope}[xshift=2cm,yshift=-1cm]
\draw (0,0) -- (-.5,1) -- (0,2) -- (.5,1) -- cycle;
\draw[black!70] (-.5,1) circle(2mm);
\draw[black!70] (.5,1) circle(2mm);
\end{scope}

\begin{scope}[xshift=4cm,yshift=-1.3cm]
\draw (0,0) -- (-.6,.8) -- (-.6,1.8) -- (0,2.6) -- (.6,1.8) -- (.6,.8) -- cycle;
\draw[black!70] (0,1.3) ellipse[x radius=1cm, y radius=1.5cm];
\end{scope}
\begin{scope}[xshift=7cm,yshift=-1.9cm]
\coordinate (A) at (-1.4,.8);
\coordinate (C) at (0,1.2);
\coordinate (D) at (.7,1);
\coordinate (E) at (1.4,.8);

\draw (0,0) -- (A) -- ($(A)+(C)$) -- ($(A)+(C)+(D)$) -- ($(A)+(C)+(D)+(E)$) -- ($(C)+(D)+(E)$) -- ($(D)+(E)$) -- (E) -- cycle;
\draw (0,0) -- (C) -- ($(A)+(C)$);
\draw (C) -- ($(C)+(E)$) -- ($(C)+(D)+(E)$);
\draw (E) -- ($(C)+(E)$);
\draw ($(A)+(C)+(D)+(E)$) -- ($(A)+(D)+(E)$) -- ($(D)+(E)$);
\draw ($(A)+(D)+(E)$) -- ($(A)+(D)$) -- (A);
\draw ($(A)+(C)+(D)$) -- ($(A)+(D)$);
\draw[black!70] (-.35,1.7) ellipse[rotate=39,x radius=1.5cm,y radius=.4cm];
\end{scope}

\begin{scope}[xshift=12cm,yshift=-2.4cm]
\coordinate (A) at (-1.4,.8);
\coordinate (B) at (-.7,1);
\coordinate (C) at (0,1.2);
\coordinate (D) at (.7,1);
\coordinate (E) at (1.4,.8);
\coordinate (F) at (0,1.8);

\draw (0,0) -- (A) -- ($(A)+(B)$) -- ($(A)+(B)+(C)$) -- ($(A)+(B)+(C)+(D)$) -- ($(A)+(B)+(C)+(D)+(E)$) -- ($(D)+(E)+(F)$) -- ($(D)+(E)$) -- (E) -- cycle;
\draw (0,0) -- (F) -- ($(A)+(B)+(C)$);
\draw (F) -- ($(F)+(E)$) -- ($(F)+(D)+(E)$);
\draw (E) -- ($(F)+(E)$);
\draw ($(A)+(B)+(C)+(D)+(E)$) -- ($(A)+(B)+(D)+(E)$) -- ($(A)+(D)+(E)$) -- ($(D)+(E)$);
\draw ($(A)+(D)+(E)$) -- ($(A)+(D)$) -- (A);
\draw ($(A)+(B)+(D)+(E)$) -- ($(A)+(B)+(D)$) -- ($(A)+(D)$);
\draw ($(A)+(B)+(C)+(D)$) -- ($(A)+(B)+(D)$) -- ($(A)+(B)$);

\end{scope}

\end{tikzpicture}
\caption{\label{fig_doubling}\scriptsize A sequence of doublings at order-convex sets satisfying the conditions of Corollary \ref{cor_doubling_crosscut}.  The final poset is the quotient lattice of Figure \ref{fig_quotient}.}
\end{centering}
\end{figure}
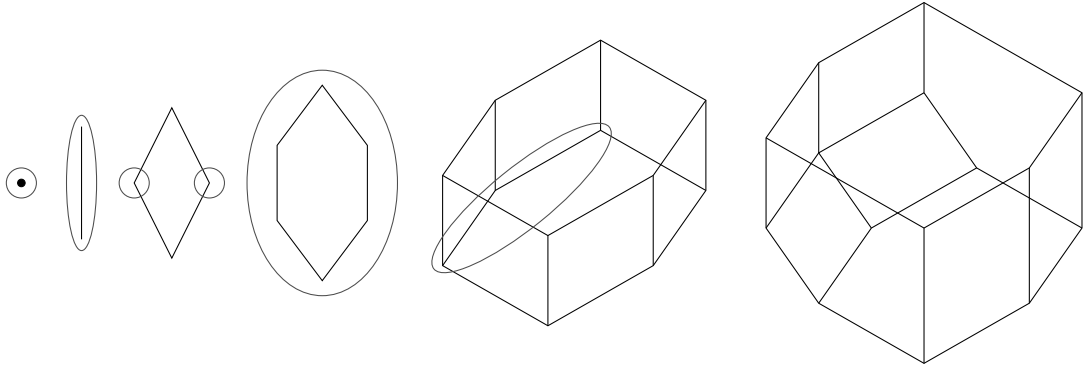

A finite lattice is distributive if and only if it may be obtained from the one-element lattice by a sequence of doublings at principal order filters.  While not every congruence-normal lattice is crosscut-simplicial, we deduce from Proposition \ref{prop_doubling_crosscut} that some doublings preserve the crosscut-simplicial property, as described in the following corollary.  An example is given in Figure \ref{fig_doubling}.

\begin{corollary}\label{cor_doubling_crosscut}
Let $C$ be an order-convex subset of a crosscut-simplicial lattice $L$.  If for $x\in C,\ y\in L-C,\ x\leq y$ the interval $[x,y]$ contains an atom not in $C$, then $L[C]$ is crosscut-simplicial.  In particular, if $C$ is an order filter, then $L[C]$ is crosscut-simplicial.
\end{corollary}

\begin{poof}
Let $I=[(x,\epsilon),(y,\epsilon^{\pr})]$ be an interval of $L[C]$.  If $I$ is an interval of type (1),(2), or (4) in Proposition \ref{prop_doubling_crosscut}, then $\Gamma(I)$ is either a simplex or its boundary.  If $I$ is of type (3), then $x\in C,\ y\in L-C,$ and $\epsilon<\epsilon^{\pr}$.  Let $A$ be the set of atoms of $[x,y]$.  By assumption, $A\cap C$ is a proper subset of $C$.  Since $\Gamma([x,y])$ is either a simplex or its boundary, the restricted complex $\Gamma([x,y])|_{A\cap C}$ is a simplex.  Therefore, $\Gamma(I)$ is isomorphic to the simplex $\{v\}*\Gamma([x,y])|_{A\cap C}$.
\end{poof}

\section{Acknowledgements}

The author thanks Tricia Hersh for explaining SB-labellings to him.  He also thanks Vic Reiner and his advisor, Pasha Pylyavskyy for their guidance.

\scriptsize
\bibliographystyle{plain}
\bibliography{bib_sb_quotients}

\end{document}